\documentclass[12pt]{amsart}

\usepackage{amsopn}
\usepackage{amsmath,amsthm,amssymb}

\newcommand{\nc}{\newcommand}

\newcommand{\fa}{{\frak a}}
\newcommand{\fD}{{\frak D}}
\newcommand{\ff}{{\frak f}}
\newcommand{\fg}{{\frak g}}
\newcommand{\fh}{{\frak h}}
\newcommand{\fH}{{\frak H}}
\newcommand{\fri}{{\frak e}}

\newcommand{\fl}{{\frak l}}

\nc{\fo}{{\frak o}}
\newcommand{\fp}{{\frak p}}
\newcommand{\fr}{{\frak r}}
\newcommand{\fs}{{\frak s}}
\newcommand{\ft}{{\frak t}}

\newcommand{\fu}{{\frak u}}
\newcommand{\fv}{{\frak v}}
\newcommand{\fw}{{\frak w}}

\newcommand{\cc}{{\Bbb C}}
\newcommand{\hh}{{\Bbb H}}
\nc{\nn}{{\Bbb N}}
\newcommand{\rr}{{\Bbb R}}
\nc{\zz}{{\Bbb Z}}

\newtheorem{teo}{Theorem}[section]
\newtheorem{lem}[teo]{Lemma}
\newtheorem{prop}[teo]{Proposition}
\newtheorem{corol}[teo]{Corollary}

\theoremstyle{remark}
\newtheorem{rem}{Remark}[section]

\nc{\lp}{\langle}
\nc{\rp}{\rangle}
\nc{\inc}{\hookrightarrow}

\title{\bf Complex structures on affine motion groups }
\author{Mar\'{\i}a L. Barberis}
\address{CIEM, FaMAF, Universidad Nacional de C\'ordoba, Ciudad Universitaria, (5000) C\'ordoba, Argentina}
\email{barberis@mate.uncor.edu}
\thanks{Both authors were partially supported 
by CONICET, ANPCyT, SECYT-UNC and ACC (Argentina).}

\author{Isabel G. Dotti}
\email{idotti@mate.uncor.edu}

\subjclass{Primary ; Secondary  }

\begin{document}
\maketitle

\begin{abstract} We study existence of complex structures on semidirect  
products $\fg \oplus _{\rho} \fv$ where $\fg$ is a real Lie algebra and 
$\rho$ is a representation of $\fg$ on $\fv$.  Our first examples,  
the Euclidean algebra $\fri(3)$ and  the Poincar\'e algebra $ \fri(2,1)$, carry complex structures obtained by deformation of  
a regular complex structure on $\fs \fl (2, \cc)$. We also exhibit a complex 
structure on the Galilean algebra $\mathcal{G}(3,1)$. We construct next a complex structure on
$\fg \oplus _{\rho} \fv$ starting with one on $\fg$ under certain compatibility 
assumptions on $\rho$.

 As an application  
 of our results we obtain that 
there exists $k\in \{ 0,1\}$ such that $(S^1)^k \times E(n)$ admits a left invariant complex structure, where $S^1$ is the circle and $E(n)$ denotes the Euclidean group. We also 
prove that the Poincar\'e group $P^{4k+3}$ has a natural left invariant complex  structure.  

In case $\dim \fg= \dim \fv$,  then 
there is an adapted complex structure on 
$\fg\oplus _{\rho} \fv$   precisely when $\rho$ determines a  flat,  torsion-free connection on $\fg$. If $\rho$ is self-dual, $\fg \oplus _{\rho}\fv$  carries a
 natural symplectic structure as well. 
If, moreover, $\rho$ comes from a metric connection 
 then  
$\fg\oplus _{\rho} \fv$ possesses a pseudo-K\"ahler structure.

We prove that the tangent bundle $TG$ of a Lie group $G$ carrying a flat torsion 
free connection $\nabla$ and a parallel complex structure possesses a  hypercomplex structure. 
More generally, by an iterative  
procedure, we can obtain Lie groups carrying a family of left invariant complex structures which generate any prescribed real Clifford algebra. 
\end{abstract}

\section{Introduction}
A complex structure on a real Lie algebra $\frak g$ is an endomorphism $J$ 
of $\frak g$ satisfying: $ J^2=-1$ and $N_J\equiv 0$, where $N_J$ denotes the Nijenhuis tensor 
corresponding to $J$ (equation \eqref{nijen} below). This integrability condition is equivalent to         a splitting of $\fg ^{\cc}=\fg \otimes _{\rr}\cc$ as 
\begin{equation}   \fg ^{\cc} = \fg^{1,0} \oplus  \fg^{0,1}            \label{split}    \end{equation}
where $\fg^{1,0}$ (resp. $\fg^{0,1}$), the $i$-eigenspace (resp. $-i$-eigenspace)  of $J$, is a complex Lie subalgebra. Observe that in case $\fg$ 
admits such a structure then its real dimension is even and 
$J$ induces a complex structure on any Lie group $G$ with Lie algebra $\fg$ 
 such that left translations by elements of $G$ are holomorphic maps on $G$.

There are two particular cases when the integrability condition is satisfied: 
\begin{itemize}
\item  $\fg^{1,0}$ and $\fg^{0,1}$ are ideals of $\fg^{\cc}$, which is equivalent to  ad$\, (x) \circ J=J \circ $ ad$\, (x)$ for all $x\in \fg$, that is, 
$(\fg, J)$ is a complex Lie algebra;
\item  $\fg^{1,0}$ and $\fg^{0,1}$ are abelian subalgebras,  which is equivalent to ad$: (\fg, J)\rightarrow (\fg \fl (\fg), R_J)$ 
being anti-holomorphic, and we say that the complex structure $J$ is abelian (cf. [\ref{bdm}], [\ref{b-d}]). 
Here, $R_J (u)= u\circ J, \; u\in \fg \fl (\fg)$.  
\end{itemize}

Samelson showed in [\ref{sam}] that any even dimensional compact Lie algebra $\fg $ admits a complex structure, which is determined by considering a root space 
decomposition of $\fg^{\cc} $.
The case of reductive $\fg$ was considered by D. Snow in [\ref{dsn}], where he classified the 
regular left invariant complex structures on $G$.
For the  solvable case there are partial results depending on $\dim _{\rr} \fg$. If $\fg$ is solvable and $\dim _{\rr}\fg =4$ the classification is complete: J. Snow 
classified in [\ref{jsn}] all four dimensinal Lie algebras $\fg$ such that $\dim 
[\fg, \fg]\leq 2$ carrying a complex structure 
 and G. Ovando extended in [\ref{ov}] this classification to the case $\dim 
[\fg, \fg] = 3$. 
If $\dim _{\rr}\fg =6$ and $\fg$ is nilpotent, S. Salamon determined 
all such Lie algebras admitting a complex structure and 
 estimated the dimension of the moduli space of complex structures on $\fg$ (cf. [\ref{ssal}]).

In the present article we focus on the problem of finding  complex structures 
on certain semidirect products $\fg \oplus _{\rho}\fv$ where $\rho: \fg \rightarrow \fg \fl(\fv)$ is a Lie algebra 
homomorphism and we look upon $\fv$ as an abelian Lie algebra. 

In Section 3
 we obtain new results concerning natural complex structures on 
automorphism groups of differential geometric structures.  
For instance, 
we exhibit left invariant complex 
structures on  the Euclidean group $E(3)$, the Poincar\'e group $P^3$ and the Galilean 
group $ G (3,1)$ and we obtain new results concerning subgroups of Aff$(\rr^n)=$GL$(\rr^n) \ltimes \rr^n$. 
Our main theorem in this section states that 
there exists $k\in \{ 0,1\}$ 
such that $(S^1)^k\times E(n)$ admits a complex structure. We point out that 
$(S^1)^k\times E(n)$ is not a complex Lie group. We  prove an analogous result for the Poincar\'e group
 $P^{4k+3}$.  

In Section 4 we restrict ourselves to representations  arising from an affine 
structure on $\fg$. Our interest 
in this situation is related to the fact that the structures  we obtain 
provide a wide variety 
of examples which deserve a better understanding due to their relation with 
symplectic and special complex geometries.  In fact, 
starting with a Lie algebra equipped with an affine structure we  obtain 
 a natural complex structure on its tangent algebra, which in addition is a special complex structure.    
Moreover, if the affine structure on $\fg$ defines a self-dual representation of $\fg$, 
the corresponding tangent algebra is special symplectic. A particular case occurs 
in the presence of a metric connection, thus obtaining a pseudo-K\"ahler structure. 
On the other hand, if $\fg$ possesses  a complex structure $J$ and an affine structure $\nabla$ such that 
$\nabla J=0$, then the tangent algebra of $(\fg, \nabla)$ admits a hypercomplex structure.
This procedure can be iterated to obtain Lie algebras possessing a family 
of complex structures generating a Clifford algebra of 
arbitrary order.

\section{Preliminaries}
We start by recalling the basic definitions.  
A  complex structure on a real Lie algebra $\fg$ is an endomorphism 
$J$\ of $\fg$ satisfying \begin{equation} J^2=-1, 
\hspace{1.5cm} N_J(x,y)=0, \; \;\; \forall x,y 
\in \fg,   \label{nijen}  \end{equation}
where $N_J$ is the Nijenhuis tensor associated to $J$:
\[    N_J(x,y)=J[x,y]-[Jx,y]-[x,Jy]-J[Jx,Jy],          \]
in other words,  
 $J$ is integrable. 
 By a hypercomplex structure we mean a 
pair of anticommuting  complex structures. More generally, a Cl$_m$-structure
 on $\fg$ is a family $\{ J_1, \dots , J_m\}$ of pairwise anticommuting 
complex structures 
on $\fg$ such that the real associative subalgebra of End$\, (\fg)$ 
generated by $\{ J_1, \dots , J_m\}$ is isomorphic to the real Clifford algebra associated to $\rr^m$ with the Euclidean inner product. Observe 
that a Cl$_2$-structure is precisely a hypercomplex structure on $\fg$. 
Cl$_m$-structures on certain solvable Lie algebras were constructed in [\ref{bdm}]. If $G$ is a Lie group with Lie algebra $\fg$ then 
a Cl$_m$-structure
 on $\fg$ can be left translated to all of $G$.  

The following lemma, whose proof is straightforward, will be used in the sequel  
when proving that $N_J\equiv 0$. 
\begin{lem} \label{lem-int} Let $\fg= \fu \oplus J\fu$ be a decomposition of $\fg$ where $\fu$ is 
a vector subspace. Then $N_J(u,v)=0, \;\; \forall \, u,v \in \fu$\ if and only if 
$N_J\equiv 0$. 
\end{lem}

We recall that an affine structure (or a left symmetric algebra structure) on a Lie algebra $\fg$   is a linear map $\rho : \fg \rightarrow 
\fg \fl (\fg)$ satisfying the following conditions: 
\begin{eqnarray}  \label{rep} \rho[x,y]&=& [\rho(x), \rho(y)], \\ 
 \label{tor-free}\rho(x)y- \rho(y)x &=& [x,y].\end{eqnarray}
These conditions are satisfied if and only if  
\[   \phi: \fg \rightarrow \fa \ff \ff(\fg), \hspace{1.2cm} x \mapsto \phi(x)=(\rho (x), x) 
 \]
is a Lie algebra homomorphism, where $\fa \ff \ff(\fg)= \fg \fl (\fg) \oplus \fg$ is the Lie 
algebra of affine motions of $\fg$ when we look upon $\fg$  as a vector space.
 If $G$ is a Lie group with Lie algebra $\fg$, then affine structures on  $\fg$ are in one to one correspondence with left invariant flat, torsion-free affine connections on $G$   
(see [\ref{bu}]). 
  
A symplectic structure on $\fg$ is a non-degenerate skew-symmetric bilinear form $\omega$ 
satisfying d$ \omega=0$, where 
\[  \text{d}\, \omega(x,y,z)= \omega (x,[y,z]) +\omega (y,[z,x]) +\omega (z,[x,y])
\] 
for $x, y, z \in \fg$.

\section{Complex structures on  affine motion Lie algebras }
Let $\fg$ be a Lie algebra and $\fv$ a $\fg$-module, that is, there exists a Lie algebra 
homomorphism    
$\rho: \fg \rightarrow \fg \fl(\fv)$. Let  $\fg \oplus _{\rho}\fv$ denote the 
 semidirect
product of $\fg$ by $\fv$, where we look upon $\fv$ as an abelian Lie algebra.  The bracket on  $\fg \oplus_{\rho} \fv$ 
is given as follows:
 \begin{equation} \label{bracket}
[(x,u),(y,v)]=([x,y],\rho(x)v -\rho(y)u) \hspace{1cm} \text{for }x, y\in \fg, \; u,v \in \fv. 
\end{equation}
We will be studying complex structures on Lie algebras of the above type. The more general 
situation $\fg \oplus_{ \rho} \fr$ where $\fr$ is a Lie algebra, $\fD (\fr)$ is the space of 
derivations of $\fr$ and $\rho : \fg  \to \fD (\fr)$ 
is a Lie algebra homomorphism   will be the object of future study. 

In case $\rho$ is a representation of $\fg$ on itself we will denote by $\fg_a$ 
the vector space underlying $\fg$, and  we refer to the semidirect product 
\begin{equation} T_{\rho} \fg := \fg \oplus _{\rho} \fg _a  \label{tang}
\end{equation}
as the \textit{tangent algebra} of $(\fg, \rho)$ (compare 
with Definition 2.9 in [\ref{dm}]).

Observe that if $ G\subset $GL$( \rr^n)$ is a Lie subgroup with Lie algebra $\fg$, then setting $\rho=$ standard 
representation on $\rr^n$, it follows that 
 the Lie algebra structure on $\fg \oplus _{\,\rho}\rr^n$  is 
 that inherited from the Lie algebra $\fa\ff\ff(\rr^n)=\fg\fl(n,\rr)\oplus \rr^n$ 
of affine motions, in other words, $\fg \oplus _{\,\rho}\rr^n$ is a Lie subalgebra of 
$\fa\ff\ff(\rr^n)$.

The next two subsections contain  some motivating examples.

\subsection{ Euclidean, Poincar\'e and Galilean algebras.}  Let $ \fs\fo (3)$ act on $\fv= \rr^3$ by the standard representation
$\rho$ via left matrix multiplication on vector columns. Then $\fs\fo(3) \oplus _{\rho}  \rr ^3\cong \fri(3)$, the  Euclidean algebra,  admits a complex structure. 
In fact, let $h=e_{12}-e_{21} , \; f_{13}=e_{13}-e_{31},  \; f_{23}=e_{23}-e_{32}$ 
be a basis of $\fs \fo(3)$, where $e_{ij}$ is the $3\times 3$ matrix with $1$ 
at the $(i,j)$ entry  and $0$ elsewhere. Let $e_i, \; 1\leq i\leq 3,$ be the 
canonical basis of $\rr^3$. Define $J$ on $\fri(3)$ as follows:
\[ Jh = e_3, \hspace{1cm} Jf_{13} = f_{23}, \hspace{1cm} Je_1=e_2.    \] 
The integrability of $J$ can be checked by direct calculation of $N_J$. This 
construction can be extended to $\fri (4k+3)$ (Theorem~\ref{isom} below).

\smallskip

Consider the Lie algebra $\fs\fo(2,1)$ of the Lorentz group acting on Minkowski 
spacetime $\rr^{2,1}$ by the standard representation $\rho$. Then $\fri (2,1)= \fs\fo(2,1)
\oplus _{\rho} \rr^{2,1}$ is the Lie algebra of the Poincar\'e group $P^3$. 
Let $s_{ij}=e_{ij} + e_{ji}, \; i\neq j$ and define a complex structure $J$ 
on $\fri (2,1)$ as follows:
\[ J h=e_{3}, \hspace{1.5cm} J s_{13}=s_{23}, \hspace{1.5cm} Je_{1}=e_{2}.   \]
In this case  the integrability of $J$ follows again by direct calculation of $N_J.$ 
This can also be extended to $\fri(4k+2,1)$ (Theorem \ref{poincare} below). 

\smallskip
It is well known that $\fs \fo (3,1) \cong \fs \fl (2, \cc)$ deforms 
into $\fri(3)$ and $\fri(2,1)$. Moreover, if $\fs \fl (2, \cc)$ has complex basis $ \{H, X_{+}, X_{-}\}$ with bracket relations
$$[H, X_+]= 2X_+,\;\; [H, X_{-}]= -2X_{-},\;\;[X_+, X_{-}]=-H $$
the complex structure given by
\[ JH = iH, \hspace{1cm} JX_{+} = iX_+, \hspace{1cm}  JX_{-} = -iX_{-}    \] 
is integrable and it is regular in the sense of [\ref{dsn}].  Its associated complex subalgebra $\fs \fo (3,1)^{1,0}$ is solvable hence the above structure  is not equivalent to the canonical one (see [\ref{ale}], where a 
classification of all homogeneous $6$-dimensional homogeneous complex manifolds is given).  Deforming 
 $$\fs \fl (2, \cc)= \fs \fo(3)^{\cc} =\fs\fo(2,1)^{\cc}$$ 
into $\fri(3)$ and $\fri (2,1)$,
it is not hard to verify that under the above deformation of $\fs \fl (2, \cc)$ the  complex structure $J$ considered above remains integrable along the deformation and yields  the complex structures considered in $\fri(3)$ and $\fri(2,1)$.

\smallskip
 Let $\mathcal{ G}(3,1)$ be the Lie algebra of the Galilean group, that is, 
\[  \mathcal{ G}(3,1)= \left\{  \;  \begin{pmatrix}  
& &  & v_1  & x_1 \\
  &A& & v_2 & x_2 \\
 & &  & v_3 & x_3 \\
0 & 0 & 0& 0& t \\
0& 0& 0& 0& 0
\end{pmatrix}   \; : \; A \in \fs \fo(3),\; x_l, v_l , t\in \rr \;      \right\}.  \] 
We show next that $\mathcal G(3,1)$ admits a complex structure. Let $\{ \; h, f_{13}, f_{23} \; \}$ 
be the 
basis of  $\fs \fo(3)$ as above  and $\text{span}\; \{    e_i, \; e_j' \; : \; 1\leq i\leq 3, 
\; 1\leq j\leq 4        \}$ a complementary subspace of  $\fs \fo(3)$ in $\fg$. Note that 
setting $v_l=0, l=1,2,3$ and $t=0$ the resulting subalgebra is $\fri (3)$.  Let $J$ be defined 
as follows:
\begin{eqnarray*} Jh &=&e_3', \hspace{1cm} Je_1'=e_2',  \hspace{1cm} Je_1=e_2,    \\ 
Jf_{13}&=&f_{23}, \hspace{.85cm} Je_3=e_4'.    \end{eqnarray*}
The integrability of $J$ follows by checking directly the vanishing of $N_J$. 
The Galilean group can be thought of as the isometry group of $\rr ^ 4 = \text{span}
\{\; e_i' \; :\; 1\leq i \leq 4 \;\}$ with the degenerate metric 
$g= dx_1^ 2 + dx_2^2 +dx_3^ 2$, and it can be obtained as a deformation of the Poincar\'e 
algebra $\fri(3,1)$ (see [\ref{gui}]).    

The preceding paragraphs can be summarized as follows:

\begin{prop}\label{euclid} The Euclidean, Poincar\'e and Galilean  algebras $\fri (3), \fri(2,1)$ 
and $\mathcal G(3,1)$ admit 
complex structures such that 
\begin{enumerate}
\item the complex structures on $\fri(3)$ and $\fri(2,1)$ are deformations 
of a regular complex structure on $\fs \fl (2, \cc)$,   
\item  the natural embedding $\fri(3)\hookrightarrow \mathcal G(3,1)$ is holomorphic.  
\end{enumerate}
\end{prop}

 We note that the above complex structure on $\fri (3)$ induces an invariant complex structure on $TS^3$, the tangent bundle of $S^3$. This follows
from the isomorphism $\fs\fo(3) \oplus _{\rho}  \rr ^3\cong \fs\fo (3) 
\oplus _{\text{ ad} }\fs\fo (3) _a$ (see \eqref{tang}), 
since the standard representation $\rho$ (via left matrix multiplication) of $ \fs\fo (3)$  on $\fv= \rr^3$  is equivalent to the adjoint representation.  
We show next that this particular example can be generalized to $TG$ for an arbitrary Lie group $G$. On the other hand, 
in [\ref{jap}] a complex structure is given to the tangent bundle of any $n$-dimensional sphere by identifying 
$T S^n$ with the 
affine hyperquadric $V_n=\{(z_1, \dots ,z_{n+1})\in C^{n+1}: \sum _i z_i^2=1\}.$

\subsection{ Tangent and cotangent bundle of a Lie group $G$ with an invariant complex structure.}      
Both, the tangent and cotangent bundle of a Lie group $G$ are Lie groups in a natural way. 
In fact, identifying the tangent bundle $TG$ with $G\times \fg$, $TG$ has a natural Lie group structure
 as the semidirect product under the adjoint representation and the corresponding 
Lie algebra is the tangent algebra of $(\fg, $ad$)$ (see \eqref{tang}).  
 Assume that  $\fg$ has a complex structure $J$ and define 
$J_+$ on $T_{\text{ad}}\fg=\fg \oplus _{\rm ad}\, \fg_a$ by $J(x,v)=(Jx,Jv)$. 
Then $J_+$ defines a complex structure on $T_{\text{ad}}\fg$. This easily follows by observing that  $N_{J_+}((x,0), (0,v))= (0, N_{J}(x,v))$, which 
vanishes by the integrability of $J$.

There is an analogous statement for the cotangent bundle $T^*G$ with the coadjoint representation. In a similar way, we get that $J_+$ defines 
a complex structure on $\fg \oplus _{{\rm ad} ^*} \fg ^*$, where $J$ on $\fg ^*$ 
is defined by $J(\alpha)(x)= -\alpha (Jx)$ 
for $\alpha \in \fg^ *, \;\; x\in \fg$. 
The integrability of $J_+$ on $T^*G$ is again a consequence of the integrability of $J$.

We summarize the preceding paragraphs as follows:

\begin{prop} \label{cor-ad} Let $G$ be a Lie group equipped with an invariant complex structure $J$. Then the tangent and cotangent bundles, $TG$ and  $T^*G$,  carry natural globally defined invariant complex structures.  
\end{prop}

\begin{rem} We note that for a complex manifold $(M,J)$,  any given torsion free connection $\nabla$ on $M$ satisfying $\nabla J=0$ allows to construct $J_+$ on $TM$ using the decomposition of $T(TM)= H\oplus V$, where $H$ is the horizontal distribution 
corresponding to the connection $\nabla$ and $V$ is the vertical distribution tangent  to the fibers. On the other hand, $T^*(M)$ always carries a natural complex structure induced from  that of $M$.
\end{rem}

\begin{rem} We show below (Proposition \ref{teo1}) that if ad is replaced by a representation $\rho$ of $ \fg$ 
on itself satisfying
certain compatibility conditions, $J_+$ as above is integrable 
on $TG=G \times \fg$ with respect to the group structure determined by $\rho$. By considering
the dual $\rho^*$ of $\rho$ we obtain the analogous result for the cotangent bundle $T^*G=
G \times \fg^*$. 
  The latter situation was considered 
in [\ref{anna}] for  $\rho$ arising from the Obata connection  associated to a hypercomplex structure on $G$. 
\end{rem}

subsection{Main results}
Let $\fv$ be a real vector space, $\dim \fv=2n$, and fix a real endomorphism 
 $I$ of $\fv$  
 satisfying $I^2=-$id.  
Given a complex structure $J$\  on $\fg$, consider on   
$\fg \oplus _{\rho}\fv$ two endomorphisms $J_+$\ and $J_-$\   
 defined by 
\begin{equation} J_{\pm}(x,v)=(Jx,\pm Iv)   \hspace{1cm}  \text{for }x \in \fg, \; v\in \fv .  \label{jaff}\end{equation}
We  study next the integrability of $J_{\pm}$ on $\fg \oplus _{\rho}\fv$. This result 
will be used throughout this section.

\begin{prop} \label{teo1} Let $(\fg, J)$ be a Lie algebra with a complex structure 
$J$ admitting a decomposition $\fg=\fg_0 \oplus \fg_1$ with $\fg_k$ $J$-stable 
subspaces, $k=0, 1$. Assume that $\fg$ acts on a complex vector space $(\fv, I)$ and the  
action $\rho: \fg \rightarrow \fg \fl(\fv)$ satisfies the following compatibility conditions:
\begin{enumerate}
\item[(i)]  $\rho(x) J = J \rho(x)$ for all $x\in \fg_0$;
\smallskip
\item[(ii)]  $\rho(Jx)Iv=\rho(x)v$ for all $x\in \fg_1, v\in \fv$. 
\end{enumerate}
Then, $J_+$ is integrable on $\fg \oplus _{\rho}\fv$. If $\fg_1=0$, then 
$J_-$ is also integrable. 
\end{prop}
\begin{proof} 
The following statements are easily proved by a direct calculation 
of $N_{J_{+}}$ and $N_{J_{-}}$, respectively:
\begin{eqnarray}
N_{J_{+}}((x,0),(0,v))&=0 \;\; \text{ if and only if } \; \;[I,\rho(x)]v
= [I,\rho(Jx)]Iv,  \\
N_{J_{-}}((x,0),(0,v))&=0 \; \; \text{ if and only if } \;\; [I,\rho(x)]v
= - [I,\rho(Jx)]Iv,
\end{eqnarray}
for $x\in \fg, \;\; v\in \fv$. The proposition now follows from these observations. 
\end{proof}

Combining the above proposition with Samelson's result (cf. [\ref{sam}]) we obtain:
\begin{corol} \label{cor-sam}
Let $\fg$ be a compact Lie algebra acting on $\cc^n$ by complex linear maps. Then 
there exists $s\in \{0, 1\}$ such that $\rr^s\oplus (\fg\oplus _{\rho}  
\cc^n)$ admits a complex structure. 
\end{corol}

\begin{rem} The above  corollary has a corresponding analogue in the 
hypercomplex case, replacing $\cc^n$ by $\hh^n$ with $\fg$ acting by 
quaternionic linear maps and allowing 
$s\in \{0, 1, 2, 3\}$ in Corollary~\ref{cor-sam}. In fact, Joyce proved in 
[\ref{joy}] 
 that given a compact Lie algebra $\fg$ there exists $s\in \{0, 1, 2, 3\}$ such that 
$\rr^s \oplus \fg$ admits a hypercomplex structure. Combining 
this fact with Proposition~\ref{teo1} we obtain the desired generalization to the hypercomplex case.   
\end{rem}
The following result is a particular case of Proposition~\ref{teo1} where $\fg=\fg_1$. 
\begin{corol}
  $\fa \ff \ff (\rr^{2n})$ admits a complex structure. \label{cor-aff}
\end{corol}
\begin{proof}Fix a real endomorphism $I$ of $\rr^{2n}$ satisfying $I^2=-$id and 
let us denote by $R_I$ 
the endomorphism of $\fg\fl(\rr^{2n})$ defined by 
\begin{equation} \label{R_I} R_I(u)=u \circ I, \hspace{2cm} u\in \fg\fl(\rr^{2n}).   
\end{equation} 
 It is straightforward that $R_I$ defines a 
complex structure on $\fg\fl(\rr^{2n})$, that is, $R_I$ is integrable (see also 
[\ref{joyc}], \S4, Example 1). In fact, identifying $\fg\fl(\rr^{2n})\cong \rr ^{4 n^2} \cong 
\cc^{2n^2}$, it turns out that 
$R_I$ is the  complex structure induced by multiplication by $i$ on $\cc^{2n^2}$. Let 
$J=R_{-I}$ and define $J_+$ on $\fa \ff \ff (\rr^{2n})=\fg\fl(\rr^{2n}) \oplus \rr^{2n}$ 
as in \eqref{jaff}. The  corollary now follows from Proposition~\ref{teo1} applied 
to $\fg=\fg\fl(\rr^{2n})$, $\fv= \rr^{2n}$ and $\rho=$ standard representation,  
by observing that  $\fg_0=0$. 

\end{proof}

We prove below the main result in this section, namely,  that the isometry group of the Euclidean $n$-space, $E(n)$, carries invariant complex structures when it is even dimensional ($n\equiv 0, 3 \; (4)$); when the dimension is odd  ($n\equiv 1, 2 \; (4)$) $S^1 \times E(n)$ carries invariant complex structures. 
  For $n\equiv 0, 1 \; (4)$  the proof follows by applying  Proposition 
\ref{teo1} with a given complex structure on 
$\fs \fo (n)$. For $n\equiv 2  \; (4)$, Proposition 
\ref{teo1} applies again  starting with 
a complex structure on $\rr z\oplus \fs\fo (n)$. 
The case $n\equiv 3  \; (4)$ differs
 from the others since $\fs \fo(n)$ is odd dimensional so it cannot be holomorphically
embedded in $\fri(n)$.  We get around this difficulty by fixing a toral subalgebra $\ft$  
in $\fs \fo(n)$, an element $h\in \ft$ and define $J$ so that $Jh$ belongs to the 
centralizer of $h$ in $\fri(n)$. 

 The Lie algebra of $E(n)$ is the Euclidean algebra $\fri (n)\cong \fs \fo (n) \oplus _{\rho}
\rr^ n$, where $\rho$ is the standard representation of $\fs \fo (n)$ on 
$\rr^n$ by left matrix multiplication on vector columns.
Before stating the theorem we introduce a complex structure $J$ on $\fs \fo(n)$, $n \equiv
0, 1    \; (4)$,   that 
will be needed later. Let $e_{ij}$ denote the square matrix with entry $1$ at the 
$(i,j)$ entry, all other entries being $0$. Then
\begin{equation} [e_{ij},e_{rs} ]= \delta _{jr}e_{is}  -\delta _{si} e_{rj}.
\end{equation}
Set $f_{ij}= e_{ij}-e_{ji}$, $h_i= f_{2i-1, 2i}$ and let $n=4k$ or $4k+1$, so that 
$\fs\fo (n)$ has rank $2k$. Define 
$J$ on $\fs \fo(n)$ such that $J^2=-1$ as follows:
\begin{equation}   Jh_i= h_{i+1}, \; 1\leq i\leq 2k, \hspace{2cm}    J f_{2j-1,l}=f_{2j,l}, \; 2\leq 2j <l\leq n.    \label{reg-J}   \end{equation}
The integrability of $J$ easily follows from the root space decomposition of $\fs\fo (n, \cc)$
 (see [\ref{hel}]). Observe that the trivial central extension $\rr z\oplus \fs \fo(4k+2)$ 
admits a complex structure which is obtained by  \eqref{reg-J} and $Jh_{2k+1}= z$.

 
\begin{teo} \label{isom} Let $\fri (n)$ be the Euclidean Lie algebra. Then there exists $s\in \{0,1\}$\ such that 
$\rr^s \oplus \fri (n)$  has a complex structure and  
the following inclusions  are holomorphic for any nonnegative integer $k$: 
\[   \fri(4k ) \hookrightarrow     \fri(4k+1)
\oplus \rr e_{4k+2}  \hookrightarrow \fri(4k+2)
\oplus \rr e_{4k+3} \hookrightarrow \fri(4k+3) . \] 
\end{teo}

\begin{proof} 
Let $n=4k$, where $k\geq 1$ is an integer and  
let $e_l$ be the vector column whose only 
non-zero entry is 
 $1$ in the  $l$th coordinate, $1\leq l \leq 4k$.  Let $I$ be the complex structure on $\rr^{4k}$ 
satisfying $I e_{2i-1} = e_{2i}, \;\; 1\leq i \leq 2k$ and $J$ the complex strucure 
on $\fs \fo(n)$ defined in equation \eqref{reg-J}. We show 
 next  that $\fg = \fs \fo(4k)$ splits as $\fg=\fg_ 0 \oplus \fg_1$
where $\fg_0$ and $\fg_1$ are $J$ invariant subspaces 
satisfying the hypothesis of 
Proposition~\ref{teo1} and therefore $J_+$ defines a complex structure 
on $\fri(4k)$. Set
\begin{gather}  
h_i= f_{2i-1, 2i},  \hspace {2cm}  1\leq i \leq 2k,   \label{cartan} \\
u_{jl}^{\pm}= f_{2j-1,2l-1}\pm f_{2j,2l},  \hspace{1.5cm}
v_{jl}^{\pm}= f_{2j-1,2l}\pm f_{2j,2l-1}, \hspace{1cm}  1 \leq j<l \leq 2k,
\label{roots}\end{gather}
and let 
\begin{gather} \fg_0 = \text{span}_{\, \rr}\, \{ \; h_i \; : \; 1\leq i\leq 2k   \;
\} 
\oplus \text{span}_{\, \rr}\, \{ \; u_{jl}^{+}, v_{jl}^{-} \; : \; 1\leq j <l \leq
2k \; \}, \\
\fg_1 = \text{span}_{\, \rr}\, \{ \; u_{jl} ^-, v_{jl}^{+} \; :\; 1\leq j <l \leq
2k \; \}. 
\end{gather}
Observe that $\fg_ 0$ is  a subalgebra isomorphic to $\fu(2k)$ and 
therefore condition (i) of the proposition is satisfied. It remains to 
show that $(J u_{jl} ^-)(Iw)= u_{jl} ^- w $ for all $w= 
\sum_{i=1}^{4k} \, w_ie_i$ in $\rr^{4k}$, $ 1\leq j<l \leq 2k$. 
 We calculate 
\begin{equation}
\begin{split}     (J u_{jl} ^-)(Iw)&= v_{jl} ^+ (\sum_{i=1}^{2k} (-w_{2i}
e_{2i-1} + w_{2i-1} e_{2i}) )   \\ &= w_{2l-1} e_{2j-1} - w_{2l} e_{2j} -    
w_{2j-1} e_{2l-1} + w_{2j} e_{2l} =  u_{jl} ^- w,           
\end{split} \end{equation}
which implies that condition (ii) of the proposition is satisfied, therefore 
$J_+(x,w)=(Jx, Iw)$ is a complex structure on $\fri(4k)$. 

We show next that the trivial central extension $\rr z \oplus \fri (4k+1)$ admits a complex structure. 
Let $ \fp= 
\text{span}_{\, \rr} \{ \; f_{i, 4k+1} \; : \; 1\leq i\leq 4k \; \}$, so that 
\[   \fri(4k+1) = \fri (4k) \oplus \fp \oplus \rr e_{4k+1}.      \]
We define $J$ on $\rr z \oplus \fri (4k+1)$ so that its restriction 
to $\fri(4k)$ is the complex structure defined above, $J$ restricted to $\fs \fo 
(4k+1)$ is given by equation \eqref{reg-J} 
and 
 \[      Je_{4k+1} =z.                            \]
 We show next that 
$J$ is integrable.   It is straightforward that $N_J(x, e_{4k+1})=0$ for $x \in 
\fri(4k)$.                
In view of Lemma~\ref{lem-int}, we only need to check that $N_J(f_{2i-1, 4k+1}, e_{2j-1})=0, 
\; \;  
1\leq i \leq k, \; 1\leq j \leq k+1$.  
It will be convenient to write down the following brackets, which will be used to carry out the calculations
\[  [f_{i, 4k+1}, e_{4k+1}]= e_i, 
\hspace{1cm}  [f_{i, 4k+1},f_{j, 4k+1}]=f_{ji}, \quad i \neq j \hspace{1cm}
[f_{i, 4k+1}, e_j]= -\delta _{ij} e_{4k+1},   
\]
for  $1\leq i, j \leq 4k$.
\begin{equation}
\begin{split} N_J(f_{2i-1, 4k+1}, e_{2j-1})  &= J[f_{2i-1, 4k+1}, e_{2j-1}] -
[Jf_{2i-1, 4k+1}, e_{2j-1}] - [f_{2i-1, 4k+1}, Je_{2j-1}] \\
& \quad - J[Jf_{2i-1, 4k+1}, Je_{2j-1}] \\
 & = -\delta _{ij} Je_{4k+1} - [f_{2i, 4k+1}, e_{2j-1}] - [f_{2i-1, 4k+1}, e_{2j}] -
J[f_{2i, 4k+1}, e_{2j}] \\
  & = -\delta _{ij} Je_{4k+1} + \delta _{ij} Je_{4k+1}=0, \hspace{2.5cm}1\leq i, j \leq 2k, \\
N_J(f_{2i-1, 4k+1}, e_{4k+1})  &= J[f_{2i-1, 4k+1}, e_{4k+1}] -
[Jf_{2i-1, 4k+1}, e_{4k+1}] - [f_{2i-1, 4k+1}, Je_{4k+1}]\\
& \quad - J[Jf_{2i-1, 4k+1}, Je_{4k+1}] \\
 & =  Je_{2i-1} - [f_{2i, 4k+1}, e_{4k+1}] - [f_{2i-1, 4k+1}, z] -
J[f_{2i, 4k+1}, z] \\
  & =  Je_{2i-1} - e_{2i}=0,   \hspace{2.5cm}1\leq i \leq 2k,      
\end{split} \end{equation}
therefore, $J$ defines a complex structure on $\rr z \oplus \fri (4k+1)$.

The case of $\fri (4k+2)$ is analogous to $\fri(4k)$. The only 
difference is that   $\fs \fo (4k+2)$ has 
odd rank equal to $2k+1$, so in this case we get that the  
trivial central extension $\rr z \oplus \fri (4k+2)$ has a complex structure. 
The 
complex structure on $\fg=\fs\fo (4k+2) \oplus \rr z$ is  that defined after \eqref{reg-J}. 
Let $\fg= \fs \fo( 4k+2) \oplus \rr z$ where $z$ is central
 and let 
$h_i, \; 1\leq i\leq 2k+1, \;  u_{jl}^{\pm}, \; v_{jl}^{\pm}, \; 
1\leq j<l \leq 2k+1 $ be defined as in \eqref{cartan} and \eqref{roots}. 
We define $I$ on $\rr ^{4k+2}$ by $I e_{2i-1} = e_{2i}, \;\; 1\leq i \leq 2k+1$. 
Proposition~\ref{teo1} applies again with 
\begin{gather}
\fg_0=\text{span}_{\, \rr}\, \{ \; h_i \; : \; 1\leq i\leq 2k +1   \; \} 
\oplus \text{span}_{\, \rr}\, \{ \; u_{jl}^{+}, v_{jl}^{-} \; : \; 1\leq j <l \leq
2k+1 \; \}
 \oplus \rr z , \\
\fg_1 = \text{span}_{\, \rr}\, \{ \; u_{jl} ^-, v_{jl}^{+} \; :\; 1\leq j <l \leq
2k +1\; \},
 \end{gather} 
  and we obtain that $J_+$   is a complex structure on $\rr z \oplus
\fri(4k+2)$.

The case of $\fri(4k+3)$ is similar to $\fri(4k+1)$. Setting $\; \fri(4k+3) = 
\fri(4k+2) \oplus \rr e_{4k+3}\oplus \fp$, where $\fp =   
\text{span}_{\, \rr} \{ \; f_{i, 4k+3} \; : \; 1\leq i\leq 4k +2\; \}\; $ it turns out
that the endomorphism $J$ on $\; \fri(4k+3)$ restricting to the complex 
structure on $\fri(4k+2) \oplus \rr e_{4k+3}$ 
defined above and satisfying $Jf_{2i-1, 4k+3}=f_{2i, 4k+3},\; 1\leq i\leq 2k+1$. The integrability of $J$ follows by  
analogous arguments to those in the case of $\fri(4k+1) \oplus \rr z$. 
\end{proof}

 \begin{rem} (i) Observe that if $(M,g)$ is a compact riemannian manifold, then the isometry 
group $I(M,g)$ is a compact Lie group, hence it always admits a complex structure if it 
is even dimensional (cf. [\ref{sam}]). 
 
(ii) If $M$ is a compact complex manifold then the group of 
holomorphic diffeomorphisms of $M$, $\fH(M)$, is a complex Lie group (cf. [\ref{boch}]). Moreover, if $M$ 
is K\"ahler-Einstein with non-zero Ricci tensor, then the Lie algebra of infinitesimal isometries 
is a real form of the Lie algebra $\fh(M)$ of the complex Lie group $\fH(M)$ (cf. [\ref{mats}]). 
 \end{rem}

The fact that the Galilean group admits a complex structure 
suggests that Theorem~\ref{isom} can possibly be extended
 to groups of transformations preserving other geometric structures. Indeed, the following theorem shows that the isometry group of Minkowski space time $\rr^{4k+2, 1}$ with 
the Lorentz metric possesses a natural complex structure. Before stating the result, we introduce
a complex structure $J$ on $\fs \fo(4k+2, 1) \oplus \rr e_{4k+3}$ which will be needed in the proof of the theorem. 
The restriction of $J$ to $\fs \fo(4k+2) \oplus \rr e_{4k+3}$
 is the complex structure defined after  \eqref{reg-J}, 
and its restriction 
to the subspace span $\{ s_{i, 4k+3}=e_{i, 4k+3}+e_{4k+3,i}, \; 1\leq i \leq 4k+2 \}$ is defined 
as follows: $$Js_{2i-1, 4k+3}= s_{2i, 4k+3}, \qquad \qquad   1\leq i \leq 2k+1.$$
The integrability of $J$ follows by observing that both, $\fs \fo(4k+3, \cc)^{1,0}$
and $\fs \fo(4k+3, \cc)^{0,1}$, are complex subalgebras of 
$\fs \fo(4k+3, \cc)$.

\begin{teo} \label{poincare} The Lie algebra $ 
\fri (4k+2,1)$  of the Poincar\'e group $P^{4k+3}$ admits a complex structure such that 
the embedding $\fri(4k+2) \oplus \rr e_{4k+3} \hookrightarrow  
\fri (4k+2,1)$ is holomorphic.
\end{teo} 
\begin{proof} Let $J$ be the complex structure on $\fs \fo(4k+2, 1) \oplus \rr e_{4k+3}$ 
defined in the previous paragraph and let $I$ on $\rr^{4k+2}$ the canonical 
complex structure defined by $Ie_{2i-1}= e_{2i}, \; 1\leq i\leq 2k+1$. It follows 
that $J_+(x,v)=(Jx,Iv), \; x \in \fs \fo(4k+2, 1) \oplus \rr e_{4k+3}, \; v\in \rr^{4k+2}$ 
is integrable on $\fri (4k+2,1)$. We apply Proposition \ref{teo1} on $\fg =\fs \fo(4k+2, 1) \oplus \rr e_{4k+3}$ as follows: 
\begin{gather*}
\fg_0=\text{span}_{\, \rr}\, \{ \; h_i \; : \; 1\leq i\leq 2k +1   \; \} 
\oplus \text{span}_{\, \rr}\, \{ \; u_{jl}^{+}, v_{jl}^{-} \; : \; 1\leq j <l \leq
2k+1 \; \}
 \oplus \rr e_{4k+3} , \\
\fg_1 = \text{span}_{\, \rr}\, \{ \; u_{jl} ^-, v_{jl}^{+} \; :\; 1\leq j <l \leq
2k +1\; \} \oplus \text{span}_{\, \rr}\, \{ \; s_{i, 4k+3} \; : \; 1\leq i \leq 4k+2 \; \}, 
 \end{gather*} 
where $h_i , u_{jl}^{\pm}, v_{jl}^{\pm}$ were defined in \eqref{cartan} and \eqref{roots}.  

The second assertion is a consequence of Theorem \ref{isom}, since the restriction
of $J_+$ to $\fri(4k+2) \oplus \rr e_{4k+3}$ is the complex structure we exhibited 
in the proof of that theorem. 
\end{proof}

\section{Complex and symplectic structures on tangent algebras }

In this section we will restrict to the case $\fv=\fg_a$, where $\fg$ is the Lie algebra of a connected Lie group $G$, and $\rho$ a representation of $\fg$ on itself.     Now, one may ask whether the endomorphism $K(x,y)=(y,-x)$ gives rise to a complex structure on $\fg \oplus_{\rho} \fg _a$. 
Related to this situation we considered in [\ref{b-d}] the following family of Lie algebras.

 Let $A$ be an associative algebra,    
 let $\fv=A_a$  and  
$\rho : A \rightarrow \fg \fl (A),\;$ 
  where $\rho (a)$ is left multiplication by $a$ in $A$.  The 
bracket on $A\oplus _{\rho} A_a$ is given as follows:
\[   [(a,b),(c,d)]= (ac-ca, ad-cb), \hspace{1cm} a,b,c,d \in A. 
 \] and a canonical  complex structure $K$ is given by 
\begin{equation} K(a,b)=(b,-a),    \hspace{1cm} a,b  \in A . \label{adapted}
 \end{equation}
 This Lie algebra was denoted by $\fa\ff \ff(A)$ in [\ref{b-d}].

Note that in the above family, the linear map $\rho$
satisfies equations \eqref{rep} and \eqref{tor-free}, in other words,  $\rho$ defines an affine structure  on the induced Lie algebra $A$ (see Section 2). In this section such representations 
will be denoted by $\nabla$ due to their relation with flat connections on Lie groups. 

\medskip

 A flat  connection on a Lie algebra $\fg$
is a linear map $\nabla : \fg \rightarrow \fg \fl (\fg), \;\; x \mapsto \nabla _x \,$ satisfying 
 \eqref{rep}, in other words, $\nabla$ is a representation. If, moreover, \eqref{tor-free} is satisfied we say that $\nabla$ is torsion-free and therefore the bracket in $ \fg$ comes 
from the left symmetric algebra structure determined by $\nabla$.

We consider next tangent algebras $T_{\nabla} \fg$  
 for some real Lie algebra $\fg$ and some flat connection $\nabla :\fg \rightarrow \fg\fl (\fg)$ (equation \eqref{tang}).

We show below that
the endomorphism $K$  given by \eqref{adapted} 
 defines a complex structure on $T_{\nabla } \fg$ precisely 
when  the flat connection $\nabla$ is 
torsion-free (compare with results in [\ref{dom}]).

  \begin{teo}  \label{teo2} Let $\fg$ be a real Lie algebra and $\nabla$ a flat affine connection 
on $\fg$. Then $K(x,y)=(y,-x)$ defines a complex structure on $T_{\nabla } \fg$ 
if and only if $\nabla$ is torsion-free. Moreover, if $\nabla$ is torsion-free 
then the connection $\nabla^ 1$ on 
$T_{\nabla } \fg$ defined by
\[        {\nabla}^1_{(x,y)}(z,w)=(\nabla_xz,\nabla_xw)
\]
is a flat, torsion-free connection such that $\nabla^1K =0$. 

Conversely, if 
$\fu$ is a Lie algebra with a complex structure $K$ 
such that 
$\fu$ decomposes as $\fu=\fg \oplus K\fg$ where 
$\fg$ is a subalgebra  and $K\fg$ is an ideal of $\fu$, then $K\fg$ 
is abelian and there 
exists a flat, torsion-free connection $\nabla$ on $\fg$  such that 
$\fu$ is the tangent algebra of $(\fg, {\nabla})$.
\end{teo}

\begin{proof}
Integrability of $K$ gives $(\nabla_xy,0)= K[(x,0)(0,y)]= ([x,y]+\nabla_yx, 0)$, that is, $\nabla$ is torsion-free. On the other hand, if $\nabla$ is torsion-free one computes:
$$[K(x,0),(y,0)] + [((x,0),K(y,0)]= (0, \nabla_yx -\nabla_xy) $$ and
$$K([(x,0),(y,0)]- [K(x,0),K(y,0)])= (0, -[x,y]),$$ hence, $K$ is integrable on $T_{\nabla} \fg$. 
 The assertion on $\nabla^1$ follows by direct calculation, 
using the properties of $\nabla$. 

To prove the converse, since $\fu=\fg \oplus K\fg$ with 
$\fg$  a subalgebra and $K\fg$  an ideal,
the map
$x\mapsto -K \circ {\rm ad}(x) \circ K$ from  $\fg: \rightarrow \fg\fl (\fg)$ 
 is a flat connection on $\fg$. Moreover,
  the integrability of $K$ toghether with the fact that $K\fg$ is an ideal gives $K\fg$ is an abelian ideal.  Finally, if $x, y \in \fg$ one has
$$ -K \circ {\rm ad}(x) \circ Ky +K\circ {\rm ad}(y) \circ Kx = [x,y]$$
since $K$ is integrable and $K\fg$ is abelian.  Thus $\nabla_x=-K\circ {\rm ad}(x)\circ K$ defines an affine structure on $\fg$ and $\fu= T_{\nabla}\fg$, as claimed.
 \end{proof}

\begin{rem}   
The above proposition is related with results in [\ref{as}] concerning 
complex product structures on Lie algebras. We recall that a complex product structure 
on $\fu$ is a pair $J, E$ where $J$ is a complex structure, $E$ is an endomorphism 
of $\fu$ anticommuting with $J$ such that $E^ 2=$id and $\fu$ splits as $\fu= \fu_+ \oplus \fu_-$ where $\fu_{\pm}$, the 
eigenspaces of $E$ of eigenvalue $\pm 1$, are subalgebras of $\fu$. It turns out that 
$\fu_-=J\fu_+$. In [\ref{as}] it is shown that when one of the eigenspaces of $E$, say $\fu_-$, is an ideal of $\fu$, then 
$\fu_-$ is abelian  so that, in view of the above proposition, $\fu$ is the tangent algebra of $(\fu_+,\nabla )$ where 
$\nabla$, the adjoint representation of $\fu$ restricted to $\fu_+$, is a flat 
torsion free connection on $\fu_+$. The above proposition says that the class of tangent algebras 
$T_{\nabla}\fg$ with $\nabla$ a flat torsion free connection on $\fg$ 
is in one to one correspondence with the class of Lie algebras $\fu$ 
admitting a complex product structure in which either $\fu_+$ or 
$\fu_-$ is an ideal of $\fu$. 
\end{rem}

\begin{corol} Let $G$ be a Lie group
admitting a left invariant flat and torsion-free connection.  Then, its tangent bundle $TG= G\times \fg$ carries 
a homogeneous complex structure which is parallel with respect to a flat torsion-free connection on $TG$ and such that the embedding of $G$ in $G\times \fg$ is totally real.
\end{corol}

As a special case of the above corollary we can quote the following various examples of Lie algebras carrying  affine structures: $k-$step nilpotent Lie algebras, $k<4$, 
$\fg\fl(\rr^{n})$, $\fa \ff \ff (\rr^n)$, Lie algebras admitting non singular derivations.
\begin{rem}
The Lie groups with Lie algebra $T_{\nabla } \fg$, $\nabla$ a flat,torsion-free connection on $\fg$, are particular examples of special complex manifolds, that is, complex manifolds $(M,J)$ together with a flat torsion free connection $\nabla$ such that $\nabla J$ is symmetric (cf. [\ref{acv}]). \end{rem}

Given a hypercomplex structure $\mathcal H={J_1,J_2}$ on $\fg$ there exists a unique 
torsion-free connection $\nabla^{\mathcal H}$, called the Obata connection associated with $\mathcal H$,  such that $\nabla^{\mathcal H}J_1=
\nabla^{\mathcal H}J_2=0$   
(see [\ref{ob}]). Furthermore, when $\nabla^{\mathcal H}$ is flat, it turns out that any  
Lie group with Lie algebra $\fg$ admits an atlas of charts such that coordinate 
changes are affine maps with quaternionic linear part (see [\ref{som}]).   
 The above theorem together with results in the previous section give:

\begin{corol} \label{cor-hyperc} Let $\fg$ be a real Lie algebra with  a flat, torsion-free connection $\nabla$ and a parallel  complex structure $J$.  
 Then $T_{\nabla } \fg$ carries a hypercomplex structure whose Obata connection
 is $\nabla^1$ as in Theorem \ref{teo2}.  
 \end{corol}
\begin{proof} The hypercomplex structure is given by 

\[    J_-(x,y)=(Jx,-Jy), \hspace{1cm} K(x,y)=(y,-x). 
\]
Observe that the condition $\nabla J=0$ is equivalent to condition (i) of Proposition \ref{teo1} 
relative to $I=J$, which ensures the integrability of $J$. The integrability of $K$ follows 
from Theorem~\ref{teo2}.  Also, $\nabla^1J_=\nabla^1K=0$, hence, by uniqueness, $\nabla^1$
 coincides with the Obata connection, as claimed. 

\end{proof}

 Combining  previous results with the above corollary we obtain a hypercomplex 
structure on $T$Aff$(\cc^{n})$, the tangent bundle of  Aff$(\cc^{n})$, and 
$T$ GL$( \rr ^{2n})$, the tangent bundle of GL$( \rr ^{2n})$. 

\begin{corol} $T$Aff$(\cc^{n})$ and $T$ GL$( \rr^{2n})$ carry natural hypercomplex structures. 
\end{corol}
\begin{proof} Consider the Lie algebra $\fa \ff \ff (\cc^n)=
\fg \fl ( \cc ^n) \oplus \cc^n$ of Aff$(\cc^{n})$ with the flat torsion free connection $\nabla$ defined 
by $\nabla_{(x,v)} (y,u)= (xy, xu)$, $x, y \in \fg\fl (\cc^n), 
\; u, v\in \cc^n$, and the complex structure $J_+$ introduced 
in Corollary \ref{cor-aff} which restricts to $\fa \ff \ff (\cc^n)$. It can be checked that $\nabla J_+ =0$ and 
the first assertion follows from Corollary \ref{cor-hyperc}. 

For the second case, 
let $J=R_I$ be the
complex structure on  $\fg\fl(\rr^{2n})$ defined in \eqref{R_I} and $\nabla$ the 
standard flat torsion free connection on $\fg\fl(\rr^{2n})$ defined by 
$\nabla_{(x,z)} (y,w)= (xy, xw)$, $x, y, z, w \in \fg\fl (\rr^{2n})$. 
We can apply again Corollary \ref{cor-hyperc} since $\nabla J=0$. 
\end{proof}

\subsection{The cotangent algebra}

Let $\nabla$ be a flat connection on $\fg$ and 
consider the contragredient representation $\nabla^*: \fg \rightarrow \fg\fl(\fg^*)$, that is:
\[    \nabla_x^*\,  \alpha=-\alpha \circ \nabla_x, \hspace{1cm} x\in \fg, \alpha \in \fg^*.       
\]

  Consider the  $2$-form $\Omega $ on the cotangent algebra $T^*_{\nabla}\fg 
:=\fg {\oplus}_{\nabla ^*}  \fg ^*$ defined by
\[  \Omega((x,\alpha ),(y,\beta))= \alpha(y)-\beta(x).
\]
 The following  result is known (compare with Theorem \ref{teo2} above); it shows that the conditions of $K$ being integrable on the tangent algebra 
and $\Omega $ being closed on the cotangent algebra 
are equivalent.

\begin{prop} [ {[\ref{Bo}]} ]   Let $\fg$ be a real Lie algebra, $\nabla$ a flat affine connection 
on $\fg$ and $\nabla^*$ the contragredient representation. Then $\Omega$  
is closed on $T^*_{\nabla}\fg$ 
if and only if $\nabla$ is torsion-free. 
\end{prop}

\medskip

We wish to study conditions which allow to induce structures from the cotangent algebra 
to   the tangent algebra. For instance, 
let $\psi: \fg \rightarrow \fg^*$ be a vector space isomorphism 
and let $\omega$ be the $2$-form on  $T_{\nabla}\fg $ induced 
by $\Omega$ and $\psi$, that is, 
\begin{equation}  \label{eqom} \omega ((x, y),(x',y')=
\Omega((x, \psi y),(x',\psi y'), \hspace{1cm} x, y, x', y'\in \fg. 
   \end{equation} 
It turns out that, in case $\psi$ is a $\fg$-module isomorphism between 
$(\fg,\nabla)$  and $(\fg, \nabla ^*)$, then the non-degenerate $2$-form $\omega$ is 
closed, that is, 
it is a symplectic structure on 
  $T_{\nabla}\fg $. In other words, the natural symplectic structure on the cotangent 
algebra can be transferred to the tangent algebra when the connection is self-dual. 

\begin{prop}\label{psi} Let $\nabla$ be a flat torsion-free connection 
on $\fg$ and assume  that $\psi:(\fg,\nabla) \rightarrow (\fg, \nabla ^*)$ is a 
$\fg$-module isomorphism. Then the $2$-form $\omega$ defined as in \eqref{eqom} is a 
symplectic structure
 on $T_{\nabla}\fg $ such that $\nabla^1 \omega=0$. 
\end{prop}

\begin{rem} Recall from [\ref{acv}] that a special symplectic structure on a manifold 
$M$ is a triple $(J, \nabla, \omega)$ where $\nabla$ is a flat torsion free connection, $J$ is a complex structure on $M$ such that $\nabla J$ is symmetric and  
$\omega$ is a parallel symplectic structure on $M$. In particular, $(J, \nabla)$ 
is a special complex structure on $M$. 
It follows from the above proposition and Theorem \ref{teo2} that $(T_{\nabla}\fg, 
K, \nabla^1,\omega)$ is a particular case of a special symplectic structure. 
\end{rem}

A particular case of the above theorem occurs in the presence of a metric connection.

Let $\fg$ be a real Lie algebra equipped with a flat, torsion-free connection 
$\nabla$ and assume that $\nabla$ is a metric connection, that is, there is a non-degenerate 
symmetric bilinear form $\langle \,\, , \, \rangle$ on $\fg$ having 
$\nabla$ as its  Levi-Civita connection. In this case there is a natural isomorphism 
between $\fg$ and its dual $\fg^*$:
\begin{equation}      \psi:\fg \rightarrow \fg^*, \hspace{1.5cm} x\mapsto x^{\flat} \label{musical}
\end{equation}  
where $x^{\flat}(y)= \langle x, y\rangle, \;\; x,y \in \fg$.  
 It follows that $\psi$ satisfies the 
hypothesis of Proposition \ref{psi} and therefore   $T_{\nabla}\fg$ carries 
a natural symplectic structure $\omega$. Observe that the complex structure $K$ from Theorem~\ref{teo2} and the symplectic form $\omega$ on 
 $T_{\nabla}\fg $ are related as follows:
\[  \langle K(x,y), (x',y') \rangle
= \omega((x,y), (x',y')), \] 
that is, $(K, \omega)$ defines a pseudo-K\"ahler structure on $T_{\nabla}\fg $. 

We summarize the above paragraph as follows: 
\begin{corol} \label{cor-music} Let $\langle \,\, , \, \rangle$ be a flat metric on $\fg$ and let $\nabla$ denote 
its Levi-Civita connection. Then $T_{\nabla}\fg$ possesses a 
natural pseudo-K\"ahler structure.  
\end{corol}

\subsection{Clifford structures}
Theorem~\ref{teo2}  suggests that one could obtain Cl$_m$-structures for arbitrary $m$ by an inductive procedure.  
In fact, this follows by observing that  if we start with a Lie algebra $\fg$ equipped 
with a flat,  torsion-free connection $\nabla$ then we can define ${\nabla}^1$ 
on $T_{\nabla } \fg$ as follows:
\[        {\nabla}^1_{(x,y)}(z,w)=(\nabla_xz,\nabla_xw)
\]
and $\nabla^1$ turns out to be a flat, torsion-free connection on 
the tangent algebra  $T_{\nabla} \fg$  such that $\nabla ^1 K =0$. Set $T^1_{\nabla} \fg=T_{\nabla} \fg$ and for $m>1$ we define inductively $ T_{\nabla}^m \fg$ to be the 
tangent algebra of $(T_{\nabla}^{m-1} \fg , \nabla ^{m-1})$, that is:
\[ 
T_{\nabla}^m \fg=T_{\nabla ^{m-1}}T_{\nabla}^{m-1} \fg, 
\hspace{.5cm}
\nabla^m_{(u,v)}(u',v')= (\nabla ^{m-1}_u u', \nabla ^{m-1}_u v').
\]
Then $\nabla^m$ is a flat, torsion-free connection on $T_{\nabla}^m \fg$.
The next theorem follows, by induction on $m$, from Theorem~\ref{teo2} and Corollary~\ref{cor-hyperc} by observing that Corollary~\ref{cor-hyperc} is the first step of the induction process. The second assertion follows from Corollary \ref{cor-music}.

\begin{teo} Let $\fg$ be a real Lie algebra carrying a flat, torsion-free 
connection $\nabla$. Then $(T_{\nabla}^m \fg, \nabla ^m)$ carries a parallel Cl$_m$-structure. 
Moreover, if $\nabla$ is a metric connection then there is a metric $g^m$ on $T_{\nabla}^m$ 
such that  $(T_{\nabla}^m, g^m, J)$ is pseudo K\"ahler for all $J \in Cl_m$ satisfying  
 $J^2=-1$. 
 
\end{teo}

\section*{References}
\begin{small}
\begin{enumerate}
\item  D.V. Alekseevsky, V. Cort\'es and C. Devchand, {\it 
Special complex manifolds}, J. Geom. Phys. {\bf 42} (1-2) (2002), 85--105; 
math.DG/9910091  
\label{acv}
\item A. Andrada and S. Salamon, {\it Complex product structures on  Lie algebras}, 
preprint 2003.    \label{as}
\item M. L. Barberis, I. G. Dotti Miatello and R. J. Miatello, 
{\it On certain locally homogeneous Clifford manifolds}, Ann. Glob. Anal. 
Geom. {\bf 13} (1995), 289--301. \label{bdm}
\item M. L. Barberis, {\it Abelian hypercomplex structures on central extensions of H-type Lie algebras}, J. Pure Appl. Algebra {\bf 158} (2001), 
15--23. \label{bar} 
\item M. L. Barberis and I. Dotti Miatello, 
{\it Hypercomplex structures on a class 
of solvable Lie groups}, Quart. J. Math. Oxford 
 (2), {\bf 47} (1996), 389--404. \label{bd} 
\item M. L. Barberis and I. Dotti, {\em Abelian complex structures on solvable Lie algebras},  math.RA/0202220, to appear in 
J. Lie theory. \label{b-d}
\item S. Bochner and D. Montgomery, {\it Groups on analytic manifolds}, Ann. of Math. 
{\bf 48} (1947), 659--669. \label{boch}
\item N. Boyom, {\it Varietes symplectiques affines}, Manuscripta Math. {\bf 64} (1989), 1--33. \label{Bo}
\item D. Burde, {\it Affine structures on nilmanifolds}, Int. J. Math. {\bf 7} (5), 1996, 
599--616.   \label{bu}
\item J. M. Dardi\'e and A. M\'edina, {\it Algebres de Lie k\"ahl\'eriennes et 
double extension}, J. Algebra {\bf 185} (3) (1996), 774--995. \label{dm}
\item P. Dombrowski, {\it On the geometry of the tangent bundle}, J. Reine Angew. Math.    
{\bf 210} (1962), 73--88 \label{dom}
\item A. Fino, {\it Cotangent bundle of hypercomplex 4-dimensional Lie 
groups}, Manuscripta Math., to appear. \label{anna}
\item Guillemin and S. Sternberg, {\it Symplectic techniques in Physics}, Cambridge University Press, 1984. \label{gui}
\item S. Helgason, {\it Differential geometry, Lie groups and 
symmetric spaces}, Academic Press, 1978. \label{hel}
\item D. Joyce, {\it Compact hypercomplex and quaternionic manifolds}, 
J. Differential Geometry {\bf 35} (1992), 743--761.  \label{joy}
\item D. Joyce, {\it Manifolds with many complex structures}, 
Quart. J. Math. Oxford 
 (2), {\bf 46} (1995), 169--184. \label{joyc}
\item Y. Matsushima {\it Sur la structure du groupe d'hom\'eomorphismes analytiques d'une 
certaine vari\'et\'e k\"ahl\'erienne}, Nagoya Math. J. {\bf 11} (1957), 145--150.  
\label{mats}
\item A. Morimoto and T. Nagano, {\it On pseudo-conformal transformation of hypersurfaces}, J. Math. Soc. Japan, {\bf 15} (1963), 289--300. \label{jap}
\item S. B. Myers and N. Steenrod, {\it The group of isometries of a riemannian manifold}, 
Ann. of Math. {\bf 40} (1939), 400--416.  \label{m-s}
\item M. Obata, {\it Affine connections on manifolds with almost complex, quaterion 
or Hermitian structures}, Japan. J. Math. {\bf 26} (1956), 43--79.  \label{ob}
\item G. Ovando, {\it Invariant complex structures on solvable real Lie groups}, Manuscripta Math. {\bf 103} (2000), 19--30.  \label{ov}
\item S. Salamon, {\it Complex structures on nilpotent Lie algebras}, J. Pure 
Appl. Algebra {\bf 157} (2001), 311--333. \label{ssal}
\item H. Samelson, {\it A class of complex analytic manifolds}, Portugal. Math.  {\bf 12} (1953), 129--132. \label{sam}
\item D. Snow, {\it Invariant complex structures on reductive Lie groups}, J. reine angew. Math. {\bf 371} (1986), 191--215. \label{dsn}
\item J. E. Snow, {\it Invariant complex structures on four dimensional solvable real Lie groups}, Manuscripta Math. {\bf 66} (1990), 397--412. \label{jsn}
\item A. Sommese, {\it Quaternionic manifolds}, Math. Ann. {\bf 212} (1975), 191--214. \label{som}
\item J. Winkelmann, {\it The classification of three-dimensional homogeneous complex manifolds}, Lecture Notes in Mathematics 1602, Springer Verlag (1995). \label{ale}
\end{enumerate}
\end{small}

\end{document}